\documentclass[10pt,reqno]{amsart}
\usepackage{amsmath,amssymb,amsfonts}
\usepackage{enumitem}
\usepackage{mathrsfs}
\newtheorem{thm}{Theorem}[section]
\newtheorem{prob}[thm]{Problem}
\newtheorem{cor}[thm]{Corollary}
\newtheorem{lem}[thm]{Lemma}
\newtheorem{prop}[thm]{Proposition}
\theoremstyle{definition}

\theoremstyle{remark}
\newtheorem{rem}[thm]{Remark}
\newtheorem{exa}[thm]{Example}

\numberwithin{equation}{section}

\newcommand{\F}{\mathbb F}
\newcommand{\Fq}{\mathbb F_q}

\newcommand{\ord}{\operatorname{ord}}
\newcommand{\rad}{\operatorname{rad}}

\newcommand{\Prim}{\mathcal P}

\begin{document}

\title[The $3$-sparsity of $X^n-1$, II]{The $3$-sparsity of $X^n-1$ over finite fields, II}
\author{Kaimin Cheng}
\address{School of Mathematical Sciences, China West Normal University, Nanchong, 637002, P. R. China}
\email{ckm20@126.com}
\subjclass[2020]{Primary 11T06}
\keywords{finite fields, cyclotomic polynomials, $3$-sparsity, irreducible factorization}
\date{\today}

\begin{abstract}
Let $q$ be a power of $2$ and let $\Fq$ be the finite field with $q$ elements.  For a positive integer $n$, the polynomial $X^n-1\in\Fq[X]$ is called $3$-sparse over $\Fq$ if every monic irreducible factor of $X^n-1$ over $\Fq$ has at most three nonzero terms.  This corrected version gives the characteristic-two classification.  Writing $n=2^\lambda m$ with $m$ odd, $X^n-1$ is $3$-sparse over $\Fq$ if and only if either $\rad(m)\mid q^2-1$, or $q=2^e$, $3\nmid e$, and $m$ lies in the exceptional $7$-family
\[
        m=7^A s_0,
        \quad A\ge1,
        \quad (s_0,7)=1,
        \quad \rad(s_0)\mid q-1,
        \quad 3\nmid s_0/\gcd(s_0,q-1),
\]
with the additional maximal $7$-adic orbit condition $\ord_{7^a}(q)=3\cdot7^{a-1}$ for $1\le a\le A$.  The latter condition is equivalent to $A=1$ or $7\nmid e$.  This condition is necessary; for example, $X^{49}-1$ is not $3$-sparse over $\F_{128}$.
\end{abstract}

\maketitle

\noindent\textbf{Correction notice.}
This version corrects the statement of the exceptional characteristic-two family involving the prime $7$.  The published version in \emph{Transactions of the London Mathematical Society} used the implication
\[
        7\nmid(q^2-1)
        \quad\Longrightarrow\quad
        \ord_{7^k}(q)=3\cdot7^{k-1},
\]
which is false for higher powers of $7$.  The correct statement includes the maximal $7$-adic Frobenius-orbit condition
\[
        \ord_{7^a}(q)=3\cdot7^{a-1}\qquad(1\le a\le A).
\]
For instance, $\ord_{49}(128)=3$, and $X^{49}-1$ is not $3$-sparse over $\F_{128}$; see Example \ref{ex:128}.

\section{Introduction}

Let $q$ be a prime power and let $\Fq$ be the finite field with $q$ elements.  The factorization of $X^n-1$ over $\Fq$ is a basic topic in finite field theory and coding theory; see, for example, \cite{Berlekamp1968,VanLint1998}.  The irreducible factors of $X^n-1$ over $\Fq$ determine the minimal cyclic codes of length $n$ over $\Fq$.  Many explicit factorization results for $X^n-1$ and related cyclotomic polynomials are known; see \cite{BlakeGaoMullin1993,BMOV2015,BRS2019,CLT2013,Meyn1996,WW2012} and the references therein.

Following Oliveira and Reis \cite{OliveiraReis2021}, we say that $X^n-1$ is \emph{$3$-sparse over $\Fq$} if every monic irreducible factor of $X^n-1$ in $\Fq[X]$ is either a binomial or a trinomial, equivalently has weight at most three.  Oliveira and Reis determined the answer for $q=2$ and $q=4$ and asked for the corresponding classification over arbitrary finite fields.

\begin{prob}\label{prob1.1}
For any prime power $q$, determine all positive integers $n$ such that $X^n-1$ is $3$-sparse over $\Fq$.
\end{prob}

\begin{prob}\label{prob1.2}
For any prime power $q$, prove or disprove that there are only finitely many primes $p$ such that $X^p-1$ is $3$-sparse over $\Fq$.
\end{prob}

The odd-characteristic case is settled in \cite{ChengOdd}: if $q$ has odd characteristic $r$ and $n=r^\lambda m$ with $(m,r)=1$, then $X^n-1$ is $3$-sparse over $\Fq$ if and only if $\rad(m)\mid q^2-1$.  The purpose of this paper is to give the corrected characteristic-two classification.

Since
\begin{equation}\label{eq:charpower}
        X^{2^\lambda m}-1=(X^m-1)^{2^\lambda},
\end{equation}
only the odd part $m$ of $n$ affects the support of the irreducible factors.  We use $\ord_m(a)$ for the multiplicative order of $a$ modulo $m$.

\begin{thm}\label{thm1.3}
Let $q=2^e$ and write
\[
        n=2^\lambda m,
        \qquad m\text{ odd}.
\]
Then $X^n-1$ is $3$-sparse over $\Fq$ if and only if one of the following two alternatives holds.
\begin{enumerate}[label=\textnormal{(\alph*)}, leftmargin=2.4em]
\item $\rad(m)\mid q^2-1$.
\item $3\nmid e$, and
\begin{equation}\label{eq:exceptional-shape}
        m=7^A s_0,
        \qquad
        A\ge1,
        \qquad
        (s_0,7)=1,
        \qquad
        \rad(s_0)\mid q-1,
\end{equation}
with
\begin{equation}\label{eq:noexcess3}
        3\nmid s_0/\gcd(s_0,q-1),
\end{equation}
and with the maximal $7$-adic orbit condition
\begin{equation}\label{eq:maximal7}
        \ord_{7^a}(q)=3\cdot7^{a-1}
        \qquad(1\le a\le A).
\end{equation}
Equivalently, under $3\nmid e$, condition \eqref{eq:maximal7} is $A=1$ or $7\nmid e$.
\end{enumerate}
\end{thm}

\begin{cor}\label{cor1.4}
For any even prime power $q$, the polynomial $X^p-1$ is $3$-sparse over $\Fq$ for only finitely many primes $p$.  More precisely, if $q=2^e$ and $S_q$ denotes the set of odd primes $p$ for which $X^p-1$ is $3$-sparse over $\Fq$, then
\[
S_q=
\begin{cases}
P(q^2-1),&\text{if }3\mid e,\\
P(q^2-1)\cup\{7\},&\text{if }3\nmid e,
\end{cases}
\]
where $P(q^2-1)$ is the set of prime divisors of $q^2-1$.
\end{cor}

\begin{proof}
Apply Theorem \ref{thm1.3} with $m=p$.  In the exceptional case one necessarily has $p=7$, and the maximal $7$-adic condition is automatic for $A=1$.
\end{proof}

\section{Preliminaries}

For $d$ coprime to $q$, let $\Phi_d(X)$ be the $d$-th cyclotomic polynomial.  We use the following standard facts; see \cite{LidlNiederreiter1997}.

\begin{lem}\label{lem:cyclo-identities}
Let $r$ be a prime and let $M$ be a positive integer with $r\nmid M$.  Then
\[
        \Phi_{Mr^a}(X)=\Phi_{Mr}(X^{r^{a-1}})\quad(a\ge1),
        \qquad
        \Phi_{Mr}(X)=\frac{\Phi_M(X^r)}{\Phi_M(X)}.
\]
\end{lem}

\begin{lem}\label{lem:degree}
If $(d,q)=1$, then every irreducible factor of $\Phi_d(X)$ over $\Fq$ has degree $\ord_d(q)$.
\end{lem}

\begin{lem}\label{lem:extension}
Let $f$ be irreducible over $\F_q$ of degree $N$.  Over $\F_{q^k}$ it factors into $\gcd(N,k)$ irreducible polynomials, each of degree $N/\gcd(N,k)$.
\end{lem}

\begin{lem}\label{lem:binomial-divides}
Let $a\in\Fq^*$ have order $M$.  Then $X^m-a$ divides $X^n-1$ over $\Fq$ if and only if $mM\mid n$.
\end{lem}

\begin{lem}\label{lem:inflation}
If $c$ and $h$ are positive integers with $\rad(h)\mid c$, then
\[
        \Phi_{ch}(X)=\Phi_c(X^h).
\]
\end{lem}

\begin{proof}
This follows by applying the identity $\Phi_{Np}(X)=\Phi_N(X^p)$ whenever $p\mid N$, one prime factor of $h$ at a time.
\end{proof}

\begin{lem}\label{lem:q2principal}
Let $(d,q)=1$ and suppose that $\rad(d)\mid q^2-1$.  Put
\[
        c=\gcd(d,q^2-1),
        \qquad
        h=d/c.
\]
Then $\ord_c(q)\in\{1,2\}$ and
\[
        \ord_d(q)=h\,\ord_c(q).
\]
Consequently $\Phi_d(X)$ factors over $\Fq$ into irreducible binomials and trinomials.
\end{lem}

\begin{proof}
Since $c\mid q^2-1$, $\ord_c(q)$ is $1$ or $2$.  The equality $\ord_d(q)=h\ord_c(q)$ follows from the usual lifting of orders at each prime power above its exponent in $q^2-1$.  Also $\rad(h)\mid c$, so Lemma \ref{lem:inflation} gives $\Phi_d(X)=\Phi_c(X^h)$.  The polynomial $\Phi_c(Y)$ factors over $\Fq$ into linear factors if $c\mid q-1$ and into quadratic factors otherwise, since all $c$-th roots lie in $\F_{q^2}$.  Substituting $Y=X^h$ gives binomials or trinomials.  Their degrees equal $\ord_d(q)$, so they are irreducible by Lemma \ref{lem:degree}.
\end{proof}

\begin{lem}\label{lem:order7}
Let $q=2^e$ and suppose $3\nmid e$.  Then
\begin{equation}\label{eq:order7}
        \ord_{7^a}(q)
        =\frac{3\cdot7^{a-1}}{\gcd(3\cdot7^{a-1},e)}
        =3\cdot7^{\max\{a-1-\nu_7(e),0\}}
        \qquad(a\ge1).
\end{equation}
In particular, \eqref{eq:maximal7} holds for all $1\le a\le A$ if and only if $A=1$ or $7\nmid e$.
\end{lem}

\begin{proof}
Since $\ord_{7^a}(2)=3\cdot7^{a-1}$ and $\ord_m(x^e)=\ord_m(x)/\gcd(\ord_m(x),e)$, formula \eqref{eq:order7} follows.  The final assertion is immediate.
\end{proof}

\begin{lem}\label{lem:phi7c}
Let $q=2^e$ with $3\nmid e$.  Let $c\mid q-1$, and let $\Prim_c$ be the set of primitive $c$-th roots of unity in $\Fq$, with $\Prim_1=\{1\}$.  Then
\begin{equation}\label{eq:phi7c}
        \Phi_{7c}(X)
        =
        \prod_{\alpha\in\Prim_c}
        \left(X^3+\alpha^2X+\alpha^3\right)
        \left(X^3+\alpha X^2+\alpha^3\right),
\end{equation}
and the displayed factors are irreducible over $\Fq$.
\end{lem}

\begin{proof}
The factorization
\[
        \Phi_7(X)=(X^3+X+1)(X^3+X^2+1)
\]
holds over $\F_2$.  Because $3\nmid e$, the two cubics remain irreducible over $\Fq$.  If $\theta$ is a root of $X^3+X+1$, then $\theta$ is a primitive seventh root of unity and $\theta^{-1}$ is a root of $X^3+X^2+1$.  For $\alpha\in\Prim_c$, the elements $\alpha\theta$ and $\alpha\theta^{-1}$ have order $7c$ and have the displayed minimal polynomials.  The product has degree $6\varphi(c)=\varphi(7c)$ and all roots are primitive $7c$-th roots.
\end{proof}

\begin{lem}\label{lem:exceptional-sufficient}
Let $q=2^e$ with $3\nmid e$.  Let
\[
        d=7^a s,
        \qquad a\ge1,
        \qquad (s,7)=1,
        \qquad \rad(s)\mid q-1.
\]
Put $c=\gcd(s,q-1)$ and $h=s/c$.  Assume $3\nmid h$ and $\ord_{7^a}(q)=3\cdot7^{a-1}$.  Then
\begin{equation}\label{eq:exceptional-factorization}
        \Phi_d(X)=
        \prod_{\alpha\in\Prim_c}
        \left(X^{3L}+\alpha^2X^L+\alpha^3\right)
        \left(X^{3L}+\alpha X^{2L}+\alpha^3\right),
        \qquad L=7^{a-1}h,
\end{equation}
and the displayed factors are irreducible trinomials over $\Fq$.
\end{lem}

\begin{proof}
Since $\rad(h)\mid c$, Lemma \ref{lem:inflation} gives $\Phi_d(X)=\Phi_{7^a c}(X^h)$.  The maximal order assumption gives $\Phi_{7^a c}(X)=\Phi_{7c}(X^{7^{a-1}})$.  Applying Lemma \ref{lem:phi7c} and then substituting $X^h$ gives \eqref{eq:exceptional-factorization}.  The order of $q$ modulo $s$ is $h$, while the order modulo $7^a$ is $3\cdot7^{a-1}$.  Because $3\nmid h$ and $(h,7)=1$, their least common multiple is $3\cdot7^{a-1}h$, which is the degree of each displayed factor.  Lemma \ref{lem:degree} gives irreducibility.
\end{proof}

\section{Necessity}

We first recall the argument which restricts possible prime divisors outside $q^2-1$.

\begin{prop}\label{prop:prime-obstruction}
Let $q=2^e$, and let $p>3$ be a prime with $p\nmid q^2-1$.  If $X^p-1$ is $3$-sparse over $\Fq$, then $3\nmid e$ and $p=7$.
\end{prop}

\begin{proof}
Put $t=\ord_p(q)$.  The cyclotomic polynomial $\Phi_p(X)$ has no binomial factor over $\Fq$: if $X^m-a$ divides $\Phi_p(X)$ and $a$ has order $M$ in $\Fq^*$, then Lemma \ref{lem:binomial-divides} gives $mM\mid p$, forcing $m=M=1$, a contradiction.  Hence every irreducible factor of $\Phi_p(X)$ is a trinomial.  Since the coefficient of $X$ in $\Phi_p(X)$ is nonzero, one such factor has the form
\[
        f(X)=X^t+aX+b,
        \qquad a,b\in\Fq^*.
\]
Let $\xi$ be a root of $f$, so $\xi$ is a primitive $p$-th root of unity.  Since $p>3$, $\xi^3$ is also primitive.  Let $G_3(X)$ be the characteristic polynomial of $\xi^3$ over $\Fq$.  Using the standard formula for characteristic polynomials of powers of roots of irreducible polynomials, one obtains
\begin{equation}\label{eq:cubic-power}
G_3(X)=
\begin{cases}
X^t+bX^{2t/3}+b^2X^{t/3}+a^3X+b^3,& t\equiv0\pmod3,\\
X^t+aX^{(2t+1)/3}+a^2X^{(t+2)/3}+a^3X+b^3,& t\equiv1\pmod3,\\
X^t+abX^{(t+1)/3}+a^3X+b^3,& t\equiv2\pmod3.
\end{cases}
\end{equation}
As $p\nmid q^2-1$, we have $t\ge3$.  The minimal polynomial of $\xi^3$ is an irreducible factor of $\Phi_p(X)$ and hence is a trinomial.  Formula \eqref{eq:cubic-power} forces $t=3$.  In this case the minimal polynomial is a scalar transform of
\[
        g(X)=X^3+X^2+1.
\]
If $3\mid e$, then $g$ splits over $\Fq$, contradicting irreducibility.  Hence $3\nmid e$.  Now a root of $g$ is a primitive seventh root of unity.  From the scalar transform just obtained, $\xi^3=b\theta$ for some primitive seventh root $\theta$ and some $b\in\Fq^*$.  Since $p\nmid q-1$, the order of $b$ is coprime to $p$, and also coprime to $7$.  Thus
\[
        p=|\xi^3|=|b\theta|=7|b|,
\]
so $p=7$.
\end{proof}

\begin{prop}\label{prop:qplus-obstruction}
Let $q=2^e$ with $3\nmid e$.  If $p$ is a prime divisor of $q+1$, then $X^{7p}-1$ is not $3$-sparse over $\Fq$.
\end{prop}

\begin{proof}
If $p=3$, then necessarily $e$ is odd.  Over $\F_2$ one has
\[
        \Phi_{21}(X)=
        (X^6+X^4+X^2+X+1)(X^6+X^5+X^4+X^2+1).
\]
Since $\gcd(6,e)=1$, Lemma \ref{lem:extension} shows that these two pentanomials remain irreducible over $\Fq$.

Now let $p>3$.  Then $p\ne7$, $\ord_p(q)=2$, and $\ord_7(q)=3$, so every irreducible factor of $\Phi_{7p}(X)$ over $\Fq$ has degree $6$.  Suppose all such factors had weight at most three.  The identity $\Phi_7(X^p)=\Phi_{7p}(X)\Phi_7(X)$ shows that a factor of $\Phi_{7p}$ has a nonzero coefficient of $X$, hence has the form $X^6+aX+b$.  Applying \eqref{eq:cubic-power} with $t=6$ to a root $\xi$ gives the minimal polynomial of $\xi^3$:
\[
        X^6+bX^4+b^2X^2+a^3X+b^3.
\]
This is an irreducible factor of $\Phi_{7p}(X)$ with five nonzero terms, a contradiction.
\end{proof}

\begin{prop}\label{prop:excess3}
Let $q=2^e$ with $e\equiv\pm2\pmod6$.  Then $X^{63}-1$ is not $3$-sparse over $\Fq$.
\end{prop}

\begin{proof}
Write $e=2t$ with $\gcd(t,3)=1$.  Let $w$ be a primitive element of $\F_4$.  A direct factorization over $\F_4$ gives the irreducible factor
\[
        X^3+X^2+X+w
\]
of $\Phi_{63}(X)$.  It has four nonzero terms.  Since $\Fq/\F_4$ has degree $t$ and $\gcd(t,3)=1$, Lemma \ref{lem:extension} shows that this cubic remains irreducible over $\Fq$.
\end{proof}

\begin{prop}\label{prop:nonmaximal}
Let $q=2^e$ with $3\nmid e$ and $7\mid e$.  Then $X^{49}-1$ is not $3$-sparse over $\Fq$.
\end{prop}

\begin{proof}
It is enough to prove the assertion over $\F_{128}$.  Let $\F_{128}=\F_2(\tau)$ with $\tau^7+\tau+1=0$.  A direct calculation gives the following irreducible factor of $\Phi_{49}(X)$ over $\F_{128}$:
\begin{equation}\label{eq:bad-cubic}
        X^3+
        (\tau^4+\tau^3+\tau^2+\tau)X^2+
        (\tau^4+\tau)X+1.
\end{equation}
It has four nonzero terms.  If $7\mid e$ and $3\nmid e$, then $\F_{128}\subseteq\Fq$ and $[\Fq:\F_{128}]=e/7$ is coprime to $3$, so the cubic \eqref{eq:bad-cubic} remains irreducible over $\Fq$ by Lemma \ref{lem:extension}.
\end{proof}

\begin{proof}[Necessity in Theorem \ref{thm1.3}]
Assume that $X^n-1$ is $3$-sparse over $\Fq$, and write $n=2^\lambda m$ with $m$ odd.  We may replace $n$ by $m$ because of \eqref{eq:charpower}.

Let $p$ be a prime divisor of $m$.  If $p\nmid q^2-1$, then Proposition \ref{prop:prime-obstruction} gives $3\nmid e$ and $p=7$.  Hence either $\rad(m)\mid q^2-1$, or $3\nmid e$ and the only prime divisor of $m$ outside $q^2-1$ is $7$.  In the first case Theorem \ref{thm1.3}(a) holds.

Assume now that the second case occurs.  Since $3\nmid e$, one has $7\nmid q^2-1$.  Write
\[
        m=7^A s_0,
        \qquad A\ge1,
        \qquad (s_0,7)=1.
\]
Every prime divisor of $s_0$ divides $q^2-1$.  No prime divisor of $s_0$ can divide $q+1$ without dividing $q-1$, because Proposition \ref{prop:qplus-obstruction} would then show that $X^{7p}-1$ is not $3$-sparse.  Therefore $\rad(s_0)\mid q-1$.

It remains to exclude excess powers of $3$ and non-maximal powers of $7$.  If $e$ is odd, then $3\mid q+1$ and $3\nmid q-1$, so Proposition \ref{prop:qplus-obstruction} with $p=3$ excludes any factor $3$ in $s_0$.  If $e$ is even and $3\nmid e$, then $e\equiv\pm2\pmod6$ and $\nu_3(q-1)=1$; an excess factor $3$ in $s_0/\gcd(s_0,q-1)$ would force $63\mid m$, contradicting Proposition \ref{prop:excess3}.  Thus \eqref{eq:noexcess3} holds.

Finally, if $7\mid e$ and $A\ge2$, then $49\mid m$, contradicting Proposition \ref{prop:nonmaximal}.  Hence $A=1$ or $7\nmid e$, which is equivalent to \eqref{eq:maximal7} by Lemma \ref{lem:order7}.  This proves Theorem \ref{thm1.3}(b).
\end{proof}

\section{Sufficiency}

\begin{proof}[Sufficiency in Theorem \ref{thm1.3}]
Again write $n=2^\lambda m$ with $m$ odd.  It is enough to factor $X^m-1$.

First suppose that $\rad(m)\mid q^2-1$.  For every divisor $d$ of $m$, Lemma \ref{lem:q2principal} shows that $\Phi_d(X)$ factors over $\Fq$ into irreducible binomials and trinomials.  Since
\[
        X^m-1=\prod_{d\mid m}\Phi_d(X),
\]
this proves that $X^n-1$ is $3$-sparse.

Now suppose that condition (b) of Theorem \ref{thm1.3} holds.  Let $d$ be a divisor of $m$.  If $7\nmid d$, then $d\mid s_0$ and $\rad(d)\mid q-1$.  Put $c=\gcd(d,q-1)$ and $h=d/c$.  Then
\[
        \Phi_d(X)=\prod_{\alpha\in\Prim_c}(X^h-\alpha),
\]
and these factors are irreducible binomials.

If $7\mid d$, write $d=7^a s$ with $1\le a\le A$ and $s\mid s_0$.  Put $c_s=\gcd(s,q-1)$ and $h_s=s/c_s$.  The condition \eqref{eq:noexcess3} gives $3\nmid h_s$, and the maximal $7$-adic condition gives $\ord_{7^a}(q)=3\cdot7^{a-1}$.  Lemma \ref{lem:exceptional-sufficient} then factors $\Phi_d(X)$ into irreducible trinomials.  Hence every irreducible factor of $X^m-1$ has at most three nonzero terms.  By \eqref{eq:charpower}, the same is true for $X^n-1$.
\end{proof}

\section{Examples and consequences}

\begin{exa}\label{ex:q4n7}
Let $q=4$ and $n=7$.  Then $e=2$, so $3\nmid e$, and $n=7^1$ satisfies the exceptional condition.  Thus
\[
        X^7-1=(X-1)(X^3+X+1)(X^3+X^2+1)
\]
over $\F_4$, with the two cubic factors irreducible.  This example shows why, in the case $e\equiv\pm2\pmod6$, an initial factor $3$ should not be required in the exceptional family.
\end{exa}

\begin{exa}\label{ex:128}
Let $q=128=2^7$ and $n=49$.  Then the numerical shape $n=7^2$ occurs, but the maximal $7$-adic condition fails because
\[
        \ord_{49}(128)=3<21.
\]
The cubic factor \eqref{eq:bad-cubic} of $\Phi_{49}(X)$ over $\F_{128}$ has four nonzero terms, so $X^{49}-1$ is not $3$-sparse over $\F_{128}$.
\end{exa}

\begin{rem}
The displayed identity
\[
        \Phi_{7^k}(X)=
        \left(X^{3\cdot7^{k-1}}+X^{7^{k-1}}+1\right)
        \left(X^{3\cdot7^{k-1}}+X^{2\cdot7^{k-1}}+1\right)
\]
is always true as a polynomial identity in characteristic two.  The point corrected here is that it is the irreducible factorization over $\F_{2^e}$ only when $\ord_{7^k}(2^e)=3\cdot7^{k-1}$.
\end{rem}

\section*{Conflict of interest}
The author declares that there is no conflict of interest.

\section*{Acknowledgment}
This work was conducted during the author's academic visit to RICAM, Austrian Academy of Sciences.  The author thanks Professor Arne Winterhof for helpful discussions and valuable suggestions.  This research was partially supported by the China Scholarship Council Fund (Grant No. 202301010002) and the Scientific Research Innovation Team Project of China West Normal University (Grant No. KCXTD2024-7).

\end{document}